\magnification=\magstep1
\input amstex
\documentstyle{amsppt}

\font\tencyr=wncyr10 % scaled \magstep1
\font\sevencyr=wncyr7 % scaled \magstep1
\font\fivecyr=wncyr5 % scaled \magstep1
\newfam\cyrfam \textfont\cyrfam=\tencyr \scriptfont\cyrfam=\sevencyr
\scriptscriptfont\cyrfam=\fivecyr
\define\hexfam#1{\ifcase\number#1
  0\or 1\or 2\or 3\or 4\or 5\or 6\or 7 \or
  8\or 9\or A\or B\or C\or D\or E\or F\fi}
\mathchardef\Sha="0\hexfam\cyrfam 58

\define\defeq{\overset{\text{def}}\to=}
\define\ab{\operatorname{ab}}
\define\pr{\operatorname{pr}}
\define\Gal{\operatorname{Gal}}

\def \isom {\overset \sim \to \rightarrow}

\define\Pic{\operatorname{Pic}}

\define\tame{\operatorname{tame}}

\define\Spec{\operatorname{Spec}}
\define\id{\operatorname{id}}

\define\Ker{\operatorname{Ker}}

\def \res{\operatorname {res}}

\def \Sec{\operatorname {Sec}}
\def\Sp{\operatorname {Sp}}

\def\char{\operatorname{char}}

\def\ur{\operatorname{ur}}
\def\adm{\operatorname{adm}}
\def\cusps{\operatorname{cusps}}
\def\cl{\operatorname{cl}}

\def\nd{\operatorname{nd}}

\def\tr{\operatorname{tr}}
\def\deg{\operatorname{deg}}
\define\Se{\frak{Sel}}
 \def \Sel{\operatorname {Sel}}

\NoRunningHeads
\NoBlackBoxes
\topmatter

\title
On the section conjecture over function fields and 
finitely generated fields
\endtitle

\author
Mohamed Sa\"\i di
\endauthor

\abstract We investigate sections of arithmetic fundamental groups of hyperbolic curves over function fields.  
As a consequence we prove that the anabelian section conjecture of Grothendieck holds over all finitely generated fields over $\Bbb Q$ if it holds 
over all number fields, under the condition of finiteness (of the $\ell$-primary parts) of certain Shafarevich-Tate groups. We also prove that if the section conjecture holds over all number fields then it holds over all finitely generated fields for curves which are defined over a number field.
\endabstract
\toc

\subhead
\S 0. Introduction
\endsubhead

\subhead
\S1 Geometrically abelian fundamental groups 
\endsubhead

\subhead
\S 2.  Sections of arithmetic fundamental groups of curves over local fields of equal characteristic $0$ 
\endsubhead

\subhead
\S 3.  Sections of arithmetic fundamental groups of curves over function fields in characteristic $0$ 
\endsubhead

\subhead
\S 4.  Proofs of Theorems A and B
\endsubhead

\subhead
\S 5.  Sections of arithmetic fundamental groups of curves over finitely generated fields in characteristic $0$ 
\endsubhead

\endtoc

\endtopmatter

\document

\footnote{Communicated by S. Mochizuki. Received January 14, 2016. Revised March 14, 2016.}
\footnote{Exeter university, e-mail address: M.Saidi\@exeter.ac.uk}

\subhead
\S 0. Introduction
\endsubhead
Let $k$ be a characteristic $0$ field and $X$ a smooth, {\it projective}, and geometrically connected {\it hyperbolic} curve
(i.e., $g(X)\ge 2$) over $k$.
Let $\pi_1(X)$ be the {\it arithmetic \'etale fundamental} group of $X$ which sits in the following exact sequence
$$1\to \pi_1(\overline X)\to \pi_1(X) @>{\pr}>> G_k\defeq \Gal (\overline k/k)\to 1,$$
where $\overline k$ is an algebraic closure of $k$ and $\overline X=X\times _k\overline k$.
In this paper we investigate continuous group-theoretic sections (i.e., splittings) $s:G_k\to \pi_1(X)$ of the projection
$\pr:\pi_1(X)\twoheadrightarrow G_k$, which
we will refer to as {\it sections} of $\pi_1(X)$.

Sections of $\pi_1(X)$ arise naturally from $k$-rational points of $X$. 
More precisely, a {\it rational} point $x\in X(k)$ determines a decomposition subgroup
$D_x\subset \pi_1(X)$, which is defined modulo conjugation by the elements of $\pi_1(\overline X)$, and which maps isomorphically to $G_k$ 
via the projection $\pr:\pi_1(X) \twoheadrightarrow G_k$. (Thus, $D_x$ determines a splitting of the above exact sequence.) 
We will refer to such a section of $\pi_1(X)$ as {\it point-theoretic}, and say that it arises from 
the rational point $x\in X(k)$.
We have a {\it set-theoretic} map
$$\varphi_X: X(k)\to  \overline {\Sec} _{\pi _1(X)},\ \ x\mapsto \varphi_X (x)=[s_x],$$
where $\overline {\Sec} _{\pi_1(X)}$ is the set of {\it conjugacy classes} of sections 
of $\pi _1(X)$, modulo conjugation by the elements of $\pi _1(\overline X)$,
and $[s_x]$ denotes the image (i.e., conjugacy class) of a section $s_x$ associated to $x\in X(k)$.

\definition {Definition 0.1} {\bf (i)}\ We say that the {\bf SC} (section conjecture) {\it holds} for $X$, over $k$, 
if the above map $\varphi_X: X(k)\to  \overline {\Sec} _{\pi _1(X)}$ is {\it bijective}.

\noindent
{\bf (ii)}\ We say that the {\bf SC} 
{\it holds} over $k$ if the {\bf SC} holds for {\it every} smooth, projective, and geometrically connected hyperbolic curve $X$ over $k$ (cf. (i)). 
\enddefinition

In his seminal letter to Faltings, Grothendieck formulated the following conjecture (cf. [Grothendieck]).

\definition {Grothendieck's Anabelian Section Conjecture (GASC)}  {\it Assume that $k$ is finitely generated over the prime field
$\Bbb Q$. Then the {\bf SC} holds over $k$}.
\enddefinition

The injectivity of the map $\varphi_X$ if $k$ is finitely generated over $\Bbb Q$, or more generally if $k$ is a sub-$p$-adic field, is well-known (cf. [Mochizuki], Theorem C). 
The statement of the {\bf GASC} is thus equivalent to the
surjectivity of the map $\varphi_X$, i.e., that {\it every section of $\pi_1(X)$ is point-theoretic} under the above assumptions 
on the field $k$. The {\bf GASC}, even over number fields, is still wide open.
More generally, one can ask:  {\it what are all fields (of characteristic $0$) for which the {\bf SC} holds}? 
In this paper we investigate {\it the section conjecture over function fields of curves in characteristic $0$}.

Given a separated, {\it smooth}, and connected {\it curve} $C$ over $k$ with function field $K\defeq k(C)$, 
an {\it abelian scheme} $\Cal A\to C$ with generic fibre
$A\defeq \Cal A\times _C\Spec K$, define the {\it Shafarevich-Tate group}
$$\Sha (\Cal A)=\Sha (A,C)\defeq \Ker \lgroup H^1(G_{K},A) \to \prod _{c\in C^{\cl}}H^1(G_{K_c},A_c)\rgroup,$$ 
where $c\in C^{\cl}$ is a closed point, $K_c$ is the completion of $K$ at $c$, $A_c\defeq A\times_{\Spec K}\Spec K_c$,
and the product is over {\it all} closed points of $C$.

\definition {Definition 0.2} Let $k'$ be a field with $\char(k')=0$. Consider the following conditions.

\noindent
{\bf (i)}\ The $\ell'$-cyclotomic character $\chi_{\ell'}:G_{k'}\to \Bbb Z_{\ell'}^{\times}$ is {\it non-Tate} in the sense of 
[Cadoret-Tamagawa] ($\S2$, Definition), meaning that the 1 dimensional $\ell'$-adic representation $\Bbb Z_{\ell'}(1)$ doesn't appear as a 
sub-representation of the representation arising from the $\ell'$-adic Tate module of an abelian variety over $k'$; $\forall$ prime integer $\ell'$.
%(a)}\ $k'$ is $\ell$-{\it cyclotomically full} for some prime integer $\ell$, meaning that the $\ell$-cyclotomic character $\chi_{\ell}:G_{k'}\to \Bbb Z_{\ell}^{\times}$ has an open image.

\noindent
{\bf (ii)}\ The {\bf SC} holds over $k'$.

\noindent
{\bf (iii)}\ Given a separated, smooth, and connected curve $C$ over $k'$, with function field $K\defeq k'(C)$, an abelian scheme $\Cal A\to C$, 
then $T\Sha(\Cal A)=0$ (cf. Notations).

\noindent
{\bf (iv)}\ Given an abelian variety $A$ over $k'$, the group of $k'$-rational points $A(k')$, and a quotient $A(k')\twoheadrightarrow D$, the following hold. 

\noindent
{\bf (a)}\ The natural map $D\to D^{\wedge}$ (cf. Notations) is {\it injective}. 

\noindent
{\bf (b)}\ The torsion group $D[N]$ is {\it finite} $\forall N\ge 1$, and $TD=0$ (cf. Notations).

\noindent
{\bf (v)}\ Given a function field $K=k'(C)$ as in (iii), $K$ admits a structure of Hausdorff topological field so that $X(K)$
becomes {\it compact} for any proper, smooth, and geometrically connected {\it hyperbolic} curve over $K$.

\noindent
{\bf (vi)}\ Given a separated, smooth, and connected curve $C$ over $k'$ with function field $K\defeq k'(C)$, a finite {\it \'etale} morphism
$\widetilde C\to C$, then the following hold. If $\widetilde C_c(k'(c))\neq \emptyset$, $\forall c\in C^{\cl}$ with residue field $k'(c)$ 
and where $\widetilde C_c$ is the scheme-theoretic inverse image of $c$ in $\widetilde C$, then $\widetilde C(K)\neq \emptyset$.

Given a field $k$ with $\char (k)=0$, we say that $k$ {\it strongly} satisfies one of the conditions (i), (ii), (iii), (iv), (v), and (vi)
above, if this condition is satisfied by any {\it finite} extension $k'/k$ of $k$. We say that the field $k$ {\it satisfies the condition $(\star)$} if 
$k$ {\it strongly} satisfies each of the conditions (i), (ii), (iii), (iv), (v), and (vi).
\enddefinition

Conditions (i), (iv), (v), and (vi) above are satisfied by {\it finitely generated} fields over $\Bbb Q$. In this case, (i) follows from the theory of weights, (iv)
follows from the Mordell-Weil and Lang-N\'eron Theorems, 
(v) follows (for the discrete topology) from  Mordell's conjecture: Faltings' Theorem and N\'eron's specialisation Theorem, 
and (vi) follows from the Hilbertian property (cf. Lemma 4.1.5) which holds for finitely generated fields. It is not known in general if condition
(iii) holds over such fields (cf. [Sa\"\i di-Tamagawa] for some partial results).

One of our main results in this paper are the following.

\proclaim {Theorem A} Let $k$ be a field with $\char (k)=0$  and $C$ a separated, smooth, and connected curve over $k$ with function field $K=k(C)$.
Assume that $k$ strongly satisfies the conditions $(i)$, $(ii)$, $(iv)$, $(v)$, and $(vi)$ in Definition 0.2. 
Let $\Cal X\to C$ be a flat, proper, and smooth relative curve, with generic fibre 
$X\defeq \Cal X\times _C\Spec K$ which is a geometrically connected hyperbolic curve over $K$, and $\Cal J\defeq \Pic^0_{\Cal X/C}$ its relative jacobian. 
Assume that $X(K)\neq \emptyset$, and $T\Sha(\Cal J)=0$. 
Then the {\bf SC} holds for $X$ (over $K$).
\endproclaim

\proclaim {Theorem B} Let $k$ be a field with $\char (k)=0$,  and $K=k(C)$ the function field of a separated, smooth, and connected curve $C$ over $k$.
Assume that $k$ satisfies the condition $(\star)$ (cf. Definition 0.2). Then $K$ strongly satisfies 
the condition $(ii)$. Equivalently, if $L/K$ is a finite extension then the {\bf SC} holds over $L$.
\endproclaim

In the case of {\it finitely generated fields} one obtains immediately from Theorem A, and Theorem B; respectively, the following Corollaries.

\proclaim {Corollary A} Assume that the {\bf SC} holds over all finitely generated fields over $\Bbb Q$ of transcendence degree $i\ge 0$.
Let $k$ be a field  with $\tr\deg_{\Bbb Q}k=i$ and $C$ a separated, smooth, and connected curve over $k$ with function field $K=k(C)$.
Let $\Cal X\to C$ be a flat, proper, and smooth relative curve, with generic fibre 
$X\defeq \Cal X\times _C\Spec K$ which is a geometrically connected hyperbolic curve over $K$, and $\Cal J\defeq \Pic^0_{\Cal X/C}$ its relative jacobian. 
Assume that $X(K)\neq \emptyset$, and $T\Sha(\Cal J)=0$. 
Then the {\bf SC} holds for $X$ (over $K$).
\endproclaim

\proclaim {Corollary B} Assume that the {\bf SC} holds over all number fields (i.e., 
all finite extensions of $\Bbb Q$) and that the condition $(iii)$ (in Definition 0.2) holds for any field $k'$ which is finitely generated over $\Bbb Q$. 
Then the {\bf SC} holds over all finitely generated fields over $\Bbb Q$.
\endproclaim

Recently, with Akio Tamagawa, we proved that the Shafarevich-Tate group  $\Sha (\Cal J)$ in Theorem A  
is indeed {\it finite} if the jacobian $J\defeq \Pic^0_{X/K}$ of $X$ is (isogenous to) a {\it constant} abelian variety, i.e., with the notations in Definition 0.2
the condition $(iii)$ holds for abelian schemes  $\Cal A\to C$ such that $A\defeq \Cal A\times _C\Spec K$  (is isogenous to an abelian variety which) 
descends (over $\overline K$) 
to an abelian variety defined over a finite extension of $k'$ (cf. [Sa\"\i di-Tamagawa], Theorem H).
As a consequence of this result, and Theorem A, one deduces (using an induction argument) the following (cf. Corollary 5.3).

\proclaim {Theorem C} Let $K$ be a finitely generated field over $\Bbb Q$ with algebraic closure $\overline K$ and $\overline {\Bbb Q}$ the algebraic 
closure of $\Bbb Q$ in 
$\overline K$. Assume that the {\bf SC} holds over all number fields. Then the {\bf SC} holds for every projective, smooth, and geometrically connected 
hyperbolic curve $X$ over $K$ such that $X\times _{\Spec K} 
{\Spec \overline K}$ is defined over $\overline {\Bbb Q}$. 
\endproclaim

Our method to prove Theorem A relies on a {\it local-global} approach  and
follows from a thorough investigation of sections of arithmetic fundamental groups of hyperbolic curves over {\it local fields of equal characteristic $0$}, and over 
{\it function fields of curves in characteristic $0$}.
Next, we explain the content of each section. 
In $\S 1$ we establish some basic facts on 
{\it geometrically abelian} fundamental groups and their sections. 
In $\S2$ we investigate (under the assumption that condition (i) holds) sections of arithmetic fundamental groups of curves over 
{\it local fields of equal characteristic $0$}. 
We observe that the section conjecture {\bf SC} {\it doesn't hold over  local fields 
of equal characteristic $0$} (cf. Lemma 2.1.3 and Proposition 2.3.1).  We discuss those sections which are point-theoretic in the case of stable curves (cf. Lemma 2.1.4 and 
Lemma 2.2.2). 
In $\S3$ we investigate sections of arithmetic fundamental groups of curves over {\it function fields}
(of transcendence degree $1$), and establish some of the basic techniques and facts in order to investigate their point-theoreticity via a {\it local-global} approach. 
In $\S4$ we prove Theorem A and explain how Theorem B can be derived from Theorem A. Finally, in $\S5$ we relate the results and techniques of $\S3$, and $\S4$,
to the results and framework in [Sa\"\i di-Tamagawa]. As a consequence, we prove some variants of Theorems A and B involving a condition of finiteness of 
the $\ell$-primary parts of certain Shafarevich-Tate groups (cf. Corollary 5.2 and Theorem 5.4), and deduce Theorem C.

Theorems A, B, and C concern sections of arithmetic fundamental groups of {\it proper} curves over function fields and finitely generated fields. One can prove similar results
for {\it non-cuspidal} sections of arithmetic fundamental groups of {\it affine} curves over function fields and finitely generated field, as well as for the birational version of the section conjecture (cf. [Sa\"\i di2]).

 \definition
{Acknowledgment} I would like to thank Akio Tamagawa for several discussions we had on the topic of this paper, as well as for 
his help with some technical details. 
Part of this work was done while the author was visiting the Research Institute of Mathematical Sciences (RIMS) of Kyoto university.
I would like very much to thank the members of the Institute for their hospitality.
\enddefinition

\definition {Notations} Given a scheme $Y$ over a field $L$ with algebraic closure $\overline L$ we write $Y_{\overline L}\defeq Y\times _{\Spec L}\Spec \overline L$ 
for the geometric fibre of $Y$.
Given a scheme $C$, a field $L$, $Y\to C$ and $\Spec L\to C$ morphisms of schemes, we write 
$Y_L\defeq Y\times _{C}\Spec L$.
For an algebraic group $G$ over a field  $L$ of characteristic $0$, with algebraic closure $\overline L$, 
we write $TG\defeq \underset{N\ge 1}\to{\varprojlim} G[N](\overline L)$ for the Tate module of $G$, where $G[N]\defeq \Ker (G@>[N]>>G)$
is the kernel of the homomorphism of multiplication by $N$.
For a profinite group $H$ we write $H^{\ab}$ for the maximal abelian quotient of $H$.
For an abelian group $D$ we write
$D^{\wedge}\defeq \underset{N\ge 1}\to{\varprojlim} D/ND$, where $ND\defeq\{N.a\ \vert\ a\in D\}$.
Given an integer $N\ge 1$, we write $D[N]\defeq \{b\in D\ \vert\  N.b=0\}$, and $TD\defeq \underset{N\ge 1}\to{\varprojlim} D[N]$ for the Tate module of $D$.
\enddefinition

\subhead
\S1 Geometrically abelian fundamental groups 
\endsubhead
Let $K$ be a field of characteristic $0$, $X\to \Spec K$ a proper, smooth, and geometrically connected {\it hyperbolic curve} over $K$. 
Let $\xi$ be a geometric point of $X$
with value in its generic point. Thus, $\xi$ determines an algebraic closure $\overline K$ of $K$.
Write $\overline \xi$ for the geometric point of $X_{\overline K}$
which is induced by $\xi$.  We have an exact sequence of \'etale fundamental groups
$$1\to \pi_1(X_{\overline K},\overline \xi)\to \pi_1(X,\xi) @>{\pr}>> G_K\to 1,\tag 1.1$$
where $G_K\defeq \Gal (\overline K/K)$.  Write 
$$\pi_1(X,\xi)^{(\ab)}\defeq \pi_1(X,\xi)/\Ker (\pi_1(X_{\overline K},\overline \xi)\twoheadrightarrow 
\pi_1(X_{\overline K},\overline \xi)^{\ab}).$$ 
We will refer to $\pi_1(X,\xi)^{(\ab)}$ as the {\it geometrically abelian} quotient of $\pi_1(X,\xi)$.

Assume that $X(K)\neq \emptyset$. Write $J\defeq \Pic^0_{X/K}$ for the jacobian variety of $X$, and $\iota:X\to J$ for the embedding
which maps a rational point $x_0\in X(K)$ to the zero section of $J$.
Then $\iota$ induces a commutative diagram of exact sequences
$$
\CD 
1  @>>> \pi_1(X_{\overline K}, \overline \xi)^{\ab} @>>> \pi_1 (X,\xi)^{(\ab)} @>>> G_K@>>> 1 \\
@.   @VVV     @VVV     @V{\id}VV  \\
1 @>>>  \pi_1(J_{\overline K}, \overline \xi) @>>> \pi_1 (J,\xi) @>>> G_K@>>> 1 \\
\endCD
\tag1.2
$$
where the vertical maps are isomorphisms, hence an identification of $G_K$-modules
$\pi_1(X_{\overline K}, \overline \xi)^{\ab} \isom \pi_1(J_{\overline K}, \overline \xi) \isom TJ$.
Let 
$$s:G_K\to \pi_1(X,\xi)$$
be a {\it section} of $\pi_1(X,\xi)$.
Then $s$ induces a {\it section} 
$$s^{\ab}: G_K\to \pi_1(X,\xi)^{(\ab)}$$
of the projection $\pi_1(X,\xi)^{(\ab)}\twoheadrightarrow G_K$.
The set of splittings of the upper sequence in diagram (1.2) is, up to conjugation by the elements of 
$\pi_1(X_{\overline K},\overline \xi)^{\ab}$, a torsor under the Galois cohomology group 
$H^1(G_K,\pi_1(X_{\overline K},\overline \xi)^{\ab})\isom H^1(G_K,TJ)$. 
We fix a {\it base point} of the torsor of splittings of this exact sequence to be
the splitting arising from the zero section of $J$, i.e., the splitting arising from the rational point $x_0\in X(K)$ (cf. above discussion). 
Then the above (conjugacy class of the) section $s^{\ab}$ corresponds to an element 
$$s^{\ab}\in H^1(G_K,TJ).$$
We will refer to $s^{\ab}$ as the {\it abelian portion} of the section $s$.

The Kummer exact sequences $0\to J[N]\to J@>N>> J\to 0$, $\forall N\ge 1$, induce an exact sequence; the so-called {\it Kummer exact sequence}
$$0\to J(K)^{\wedge}\to H^1(G_K,TJ)\to TH^1(G_K,J)\to 0.\tag 1.3$$
We will identify $J(K)^{\wedge}$ with its image in $H^1(G_K,TJ)$ via the above Kummer map $J(K)^{\wedge}\to H^1(G_K,TJ)$.
Note that there exist natural maps $X(K)@>{\iota}>>  J(K) \to J(K)^{\wedge}$,
where for $x\in X(K)$ the image $\iota (x)$ is the class $[x-x_0]$ of the degree $0$ divisor $x-x_0$.

\definition {Definition 1.1}
Let $\eta \in H^1(G_K,TJ)$. 

\noindent
{\bf (i)}\ We say that $\eta$ is {\it pro-geometric} if $\eta$ lies in the subgroup
$J(K)^{\wedge}$ of $H^1(G_K,TJ)$. 

\noindent
{\bf (ii)}\ We say that $\eta$ is {\it geometric} if $\eta$ is in the image of the natural composite homomorphism
$J(K)\to J(K)^{\wedge}\to H^1(G_K,TJ)$. 

\noindent
{\bf (iii)}\ We say that $\eta$ is {\it point-theoretic} if $\eta$ is in the image of the composite map
$X(K)@>{\iota}>> J(K)\to  J(K)^{\wedge}\to H^1(G_K,TJ)$.
\enddefinition

The following Lemma follows easily from the various Definitions.

\proclaim {Lemma 1.2} 
Suppose that the section $s=s_x$, $x\in X(K)$, is point-theoretic (cf. $\S0$).
Then $s^{\ab}\in H^1(G_K,TJ)$ is point-theoretic (cf. Definition 1.1(iii)), and $s^{\ab}$ is the image of $x$ via the composite map 
$X(K)@>{\iota}>>  J(K) \to J(K)^{\wedge}\to  H^1(G_K,TJ)$.  
In particular, $s^{\ab}$ is pro-geometric and geometric 
(cf. Definition 1.1(i)(ii)). 
\endproclaim

\subhead
\S 2. Sections of arithmetic fundamental groups of curves over local fields of equal characteristic $0$ 
\endsubhead
In this section $K$ is a {\it complete discrete valuation field} of {\it equal characteristics} $0$, $\Cal O_K$ its valuation ring,
and $k$ its residue field. We use the notations introduced in $\S0$ and $\S1$.
Moreover, we {\it assume that $k$ satisfies the condition (i) in Definition 0.2}, unless we specify otherwise.

Let $X\to \Spec \Cal O_K$ be a flat, proper, {\it stable}, and geometrically connected (relative) curve over $\Cal O_K$, with $X_K$ {\it smooth}. 
Assume that the irreducible components $\{X_i\}_{i\in I}$ of $X_k=\sum _{i\in I}X_i$ are {\it smooth}, {\it geometrically connected}, and the singular points 
$\{x_j\}_{j\in J}$ of $X_k$ are $k$-{\it rational}. 
Let $\xi$ (resp. $\xi'$) be a geometric point of $X$ with value in the generic  point of $X$ (resp. with value in the generic point of some irreducible component
$X_{i_0}$ of $X_k$). Thus, $\xi$ (resp. $\xi'$) determines an algebraic closure $\overline K$ (resp. $\overline k$) of $K$ (resp. $k$). 
We have natural exact sequences of {\it arithmetic \lq\lq{admissible}\rq\rq fundamental groups}
$$1\to \pi_1(X_{\overline K}, \overline \xi) \to \pi_1 (X_{K},\xi) \to G_K\to 1,$$
and
$$1\to \pi_1(X_{\overline k}, \overline \xi')^{\adm} \to \pi_1 (X_{k},\xi')^{\adm} \to G_K\to 1,$$
where $G_K\defeq \Gal (\overline K/K)$, and the geometric point $\overline \xi$
(resp. $\overline \xi'$) is naturally induced by $\xi$ (resp. $\xi'$). Here, the superscript \lq\lq{$\adm$}\rq\rq means {\it admissible fundamental group}
(cf. [Mochizuki1], Definition 2.7, for more details on the definition of $\pi_1(X_{\overline k}, \overline \xi')^{\adm}$. We also refer to [Mochizuki1], the discussion before Definition 2.1, 
as well as the discussion between  Definition 2.5 and Definition 2.7, for the precise definition of the log structures on $X_k$, $\Spec k$, $\Spec \Cal O_K$, 
$X_k\to \Spec k$, and $X\to \Spec \Cal O_X$ involved in the definition of $\pi_1(X_{\overline k}, \overline \xi')^{\adm}$). 
Moreover, after a suitable choice of $\xi$ and $\xi'$ we have a commutative diagram of exact sequences
$$
\CD
1  @>>> \pi_1(X_{\overline K}, \overline \xi) @>>> \pi_1 (X_{K},\xi) @>>> G_K@>>> 1 \\
@.   @VVV     @V{\Sp}VV     @V{\id}VV  \\
1  @>>>   \pi_1(X_{\overline k}, \overline \xi')^{\adm} @>>>  \pi_1 (X_{k},\xi')^{\adm} @>>> G_K @>>> 1  \\
\endCD
\tag 2.1$$
where the middle and left vertical maps are continuous homomorphisms of \lq\lq{\it specialisation}\rq\rq, which are isomorphisms 
since $\char(k)=0$ (cf. loc. cit. Definition 2.4, Lemma 2.6, and the discussion before Definition 2.7). 
We assume that $X_K$ is {\it hyperbolic}, i.e., $g(X_{K})\ge 2$.
We have an exact sequence
$1\to I_K\to G_K\to G_k\defeq \Gal (\overline k/k)\to 1$
where $I_K$ is the {\it inertia} group. Moreover, $I_K\isom \hat {\Bbb Z} (1)$ where the \lq\lq{$(1)$}\rq\rq is a Tate twist.

\subhead
2.1. The good reduction case
\endsubhead
Assume that $X$ is {\it smooth}. 
In this case we have a commutative diagram of exact sequences of arithmetic fundamental groups
$$
\CD 
1  @>>> \pi_1(X_{\overline K}, \bar \xi) @>>> \pi_1 (X_{K},\xi) @>>> G_K@>>> 1 \\
@.   @VVV     @V{\Sp}VV     @VVV  \\
1  @>>>   \pi_1(X_{\overline k}, \bar \xi') @>>>  \pi_1 (X_{k},\xi') @>>> G_k @>>> 1  \\
\endCD
\tag 2.2$$
where the middle vertical map is the surjective homomorphism of {\it specialisation} (defined up to conjugation, cf. [Grothendieck1], Expos\'e X, $\S2$), 
the left vertical map is an isomorphism (since $\char(k)=0$), and the right vertical map is a
surjection.

\proclaim {Lemma 2.1.1} The following hold.

\noindent
(i)\ The right square in diagram (2.1) is cartesian.

\noindent
(ii)\  The projection $\pi_1(X_{K},\xi)\twoheadrightarrow G_K$ induces
a natural isomorphism $\Ker (\Sp)\isom I_K$ (i.e., maps $\Ker (\Sp)$ isomorphically onto $I_K$).
\endproclaim

\demo {Proof} The proof of (i) is similar to the proof of Lemma 3.3.2 in [Sa\"\i di]. Assertion (ii) is clear in light of (i).
\qed
\enddemo

Let 
$s:G_K\to \pi_1(X_{K},\xi)$ 
be a {\it section} of $\pi_1 (X_{K},\xi)$ so that it induces, by composing
with the specialisation homomorphism $\Sp:   \pi_1(X_{K},\xi)\to \pi_1(X_k,\xi')$,
a continuous homomorphism $s'\defeq \Sp\circ s:G_K\to \pi_1(X_k,\xi')$.

\proclaim {Lemma 2.1.2} The equality $\Ker (s')=I_K$ holds. In particular, $s'$ factorizes through $G_k$ and
induces a section
$\tilde s:G_k\to \pi_1(X_k,\xi')$ of $\pi_1(X_k,\xi')$. 
\endproclaim

\demo {Proof} It suffices to show that the image $s'(I_K)$ of the inertia subgroup in $\pi_1(X_k,\xi')$ is trivial.
This image is contained in $\pi_1(X_{\overline k},\bar \xi')$ by diagram (2.1).  A standard (well-known) weight argument,
using the fact that $k$ satisfies the condition (i) in Definition 0.2,
shows that this image must be trivial (cf. [Hoshi-Mochizuki], Lemma 1.6). 
\qed
\enddemo

Assume that the section $s$ is {\it point-theoretic}, i.e., $s=s_x:G_K\to \pi_1(X_{K},\xi)$ is associated to a {\it rational point}
$x\in X(K)$ (cf. $\S0$). Let $\overline x\in X(k)$ be the {\it specialisation} of the point $x$. Then one verifies easily that the section $\tilde s:G_k\to \pi_1(X_k,\xi')$
of $\pi_1(X_k,\xi')$ which is induced by $s$ (cf. Lemma 2.1.2) 
is {\it point-theoretic} and arises from the $k$-rational point $\overline x$, i.e., $[\tilde s]=[s_{\overline x}]$ holds in $\overline {\Sec} _{\pi _1(X_k,\xi')}$.

\proclaim {Lemma 2.1.3} We use the same notations as above. Let $x'\in X(K)$ be a rational point which specialises in $\overline x$. Then $[s_x]=[s_{x'}]$. In particular,
the map $\varphi_X:X(K)\to \overline {\Sec}_{\pi_1(X_{K},\xi)}$ (cf. $\S0$) is not injective.
\endproclaim

\demo {Proof} Indeed, it follows immediately from Lemma 2.1.1(i) and Lemma 2.1.2 that a section  $s:G_K\to \pi_1(X_{K},\xi)$ of  
$\pi_1(X_{K},\xi)$ is uniquely determined by the continuous
homomorphism $s'\defeq \Sp\circ s:G_K\twoheadrightarrow G_K/I_K\to  \pi_1(X_k,\xi')$ that it induces. In particular, all rational points $x'\in X(K)$ which specialise in
$\overline x$ (there are infinitely many such points $x'$) give rise to the same conjugacy class of sections of $\pi_1(X_{K},\xi)$, from which the second assertion follows.
\qed
\enddemo

Conversely we have the following.

\proclaim {Lemma 2.1.4}  Assume that the section $\tilde s:G_k\to \pi_1(X_k,\xi')$ of 
$\pi_1(X_k,\xi')$ which is induced by $s$ (cf. Lemma 2.1.2) is point-theoretic, i.e., 
$\tilde s=s_{\overline x}$ for some $k$-rational point $\overline x\in X(k)$. Then the section $s$ is point-theoretic, i.e., $s=s_x$ for some (non unique) $x\in X(K)$ which 
specialises in the point $\overline x$. 
\endproclaim

\demo {Proof} Let $x\in X(K)$ be a rational point which specialises
in $\overline x$ (such a point $x$ exists since $X$ is smooth). Then $[s]=[s_x]$ holds by the same argument used in the proof of Lemma 2.1.3. 
\qed
\enddemo

\subhead
2.2. The bad reduction case
\endsubhead
In this section, and in addition to our assumptions, we will {\it assume that $k$ satisfies the condition (iv)(a) in Definition 0.2} unless we specify otherwise.
We suppose that $X_k$ is {\it singular}. Recall that the irreducible components $\{X_i\}_{i\in I}$ of $X_k=\sum _{i\in I}X_i$ are {\it smooth}, {\it geometrically connected}, and the singular points $\{x_j\}_{j\in J}$ of $X_k$ are all $k$-{\it rational}.

Let $X_i$ be an irreducible component of $X_k$ and
$D_{X_i}\subset \pi_1 (X_{k},\xi')^{\adm}$ 
a {\it decomposition group} associated to $X_i$. Thus, $D_{X_i}$ is the decomposition group of an irreducible component 
of the special fibre of the universal admissible cover $\widetilde X$ of $X$
which lies above the component $X_i$, and $D_{X_i}$ is only defined up to conjugation 
(cf. [Mochizuki1], discussion before Proposition 4.1).  Let $\overline X_i$ be the (unique, since $X_i$ is geometrically connected) 
irreducible component of $X_{\overline k}$ above $X_i$. Then we have the following commutative diagram.

$$
\CD
1@>>>   D _{\overline X_i}@>>> D_{X_i}@>>> G_K @>>> 1\\
@. @VVV    @VVV   @VVV   \\
1@>>>  \pi_1(\overline X_i, *)^{\adm} @>>>  \pi _1 (X_i,*)^{\adm}  @>>> G_K @>>> 1\\
@.  @VVV   @VVV   @VVV  \\
1@>>>  \pi_1(\overline U_i, *)^{\tame} @>>>  \pi _1 (U_i,*)^{\tame}  @>>> G_k @>>> 1 \\
\endCD
\tag 2.3$$
Here, $D_{\overline X_i}$ is defined so that the upper horizontal sequence is exact, $*$ denote suitable base points, 
and $\pi_1(X_i,*)^{\adm}$ denotes the {\it admissible} fundamental group of $X_i$ which is marked by the {\it cusps}, i.e., the double points of 
$X_k$ lying on $X_i$ (cf. [Mochizuki1], Proposition 4.2 and the discussion before it which explains the log structure on $X_i$ and the choice of the base point $*$
involved in the definition of $\pi_1(X_i,*)^{\adm}$),
 $\pi_1(\overline X_i, *)^{\adm}\defeq  \Ker (\pi _1 (X_i,*)^{\adm} \twoheadrightarrow G_K)$, $U_i\defeq X_i
\setminus \{\cusps\}$, and $\overline U_i\defeq \overline X_i\setminus \{\cusps\}$, respectively. 
The superscript \lq\lq{tame}\rq\rq means {\it tame} fundamental group.
The left and middle upper vertical maps are isomorphisms defined up to conjugation (cf. loc. cit.), 
and the lower right square is cartesian. 
Note that $\pi_1(\overline U_i, *)^{\tame}=\pi _1 (\overline U_i,*)$,
and $\pi_1(U_i, *)^{\tame}=\pi _1 (U_i,*)$, since $\char (k)=0$.

Let 
$s:G_K\to \pi_1(X_{K},\xi)$ 
be a {\it section} of $\pi_1 (X_{K},\xi)$, which induces a section
$s'\defeq \Sp\circ s:G_K\to \pi_1(X_k,\xi')^{\adm}$
of the projection $\pi_1(X_k,\xi')^{\adm}\twoheadrightarrow G_K$ (cf. diagram (2.1)). 
Suppose that $X$ is {\it regular}
and $s$ is {\it point-theoretic}, i.e., 
$s=s_{x}$ arises from a rational point $x\in X(K)$ (cf. $\S0$). 
Then the $K$-rational point $x$ {\it specialises} in a rational point $\overline x\in X(k)$ which is a {\it smooth} 
point of $X_k$ and lies on  a unique irreducible component  
$X_i$ of $X_k$ (cf. [Liu], Corollary 9.1.32). 
Moreover, it follows from the various definitions that 
$s(G_K)\subset D_{X_i}\subset \pi_1 (X_{k},\xi')^{\adm}$ 
holds, where $D_{X_i}$ is a {\it decomposition group} associated to $X_i$ (cf. above discussion).
In particular, the section $s=s_x$ induces a section 
$s_i:G_K\to \pi_1(X_i,*)^{\adm}$ 
of the projection $\pi_1(X_i,*)^{\adm}\twoheadrightarrow G_K$, and a continuous homomorphism
$\tilde s_i:G_K\to \pi _1 (U_i,*)^{\tame}$
(cf. diagram (2.3)).

\proclaim {Lemma 2.2.1} The following hold.

\noindent
(i)\  The section $s_i$ is unramified, i.e., $\tilde s_i(I_K)=\{1\}$. In particular,  
$\tilde s_i$ induces a section  
$\overline s_i:G_k\to \pi _1 (U_i,*)^{\tame}$  
of $\pi _1 (U_i,*)^{\tame}$.

\noindent
(ii)\  The section $\overline s_i:G_k\to \pi _1 (U_i,*)^{\tame}$ in (i) is point-theoretic and arises from the 
$k$-rational point $\overline x\in U_i(k)$. Moreover, the section $\overline s_i$ is non-cuspidal, i.e., 
$\overline s_i(G_k)$ is not contained in a decomposition group associated to a cusp.
\qed
\endproclaim

\demo {Proof} 
Note that $U_i$ is hyperbolic since $X$ is stable.
Assertion (i) follows from the assumptions (i) and (iv)(a). 
First, assumption (iv)(a) implies that the closed points of $X_i$ are uniquely determined by their corresponding (conjugacy classes of) decomposition groups in 
$\pi _1 (U_i,*)^{\tame}$, such a decomposition group is self-normalising in $\pi _1 (U_i,*)^{\tame}$, and no non-cuspidal decomposition group is contained 
in a cuspidal decomposition group (cf. the arguments used in the proof of Theorem 1.3 in [Mochizuki2], and [Tamagawa], Proposition 2.8(i)). 
Second, if $\tilde s_i(I_K)$ is non trivial then one shows, using the condition (i), that $\tilde s_i(I_K)\subseteq  \pi_1(\overline U_i, *)^{\tame}$ 
would be a non-trivial (necessarily torsion-free) procyclic group contained in an inertia group $I_y$ at a cusp $y\in X_i\setminus U_i$ (cf. [Hoshi-Mochizuki], Lemma 1.6).
Moreover, $\tilde s_i(G_K)$ (which normalises $\tilde s_i(I_K)$) is contained in a decomposition group associated to $y$ and the latter would contain a decomposition group associated to $\overline x$, which is a contradiction (cf. above discussion). 
Assertion (ii) follows easily.
\qed
\enddemo

Conversely, suppose that the section $s$ satisfies $s(G_K)\subset D_{X_i}$, 
i.e., the image of $s$ is contained in a decomposition group
associated to an irreducible component $X_i$ of $X_s$. 
Thus, $s$  induces a section 
$s_i:G_K\to \pi_1(X_i,*)^{\adm}$
 of the 
projection $\pi_1(X_i,*)^{\adm}\twoheadrightarrow G_K$ which induces a homomorphism
$\tilde s_i:G_K\to \pi _1 (U_i,*)^{\tame}$ 
(cf. diagram (2.2)). Assume further that $s_i$ is {\it unramified}, i.e., $\tilde s_i(I_K)=\{1\}$. Then $\tilde s_i$ 
induces naturally a section 
$\overline s_i:G_k\to \pi _1 (U_i,*)^{\tame}$ 
of $\pi _1 (U_i,*)^{\tame}$ (cf. diagram (2.2)). 

\proclaim {Lemma 2.2.2} With the same notations/assumptions as above 
suppose that the section $\overline s_i:G_k\to \pi _1 (U_i,*)^{\tame}$ is point-theoretic and arises from a rational point $\overline x\in U_i(k)$. 
Then the section $s$ is point-theoretic and
arises from a (non unique) rational point $x\in X(K)$ which specialises in the point $\overline x\in U_i(k)$.    
\endproclaim

\demo{Proof} Similar to the proof of Lemma 2.1.4.
\qed
\enddemo

\subhead 2.3
\endsubhead
Next, we provide examples of sections of arithmetic fundamental groups of hyperbolic curves
over local fields of equal characteristic $0$ which are {\it not point-theoretic}. We use the notations and assumptions in 2.2.
Let $X$ be a {\it regular and stable $\Cal O_K$-curve} as above satisfying the conditions in 2.2. 
Let $x_j\in X(k)$ be a $k$-{\it rational double point} of $X_k$ and write $D_{x_j}\subset \pi_1(X_k,\xi')^{\adm}$
for the {\it decomposition group} of $x_j$.
Thus, $D_{x_j}$ is the decomposition group of a closed point
of the special fibre of the universal admissible cover $\widetilde X$ of $X$ which lies above the double point $x_j$, and $D_{x_j}$ is only defined up to conjugation 
(cf. [Mochizuki1], $\S 5$ and $\S 6$).  We have an exact sequence
$$1 \to \Delta _{x_j}  \to D_{x_j}  \to G_K \to 1,$$
where  $\Delta _{x_j}\defeq \Ker(D_{x_j}\twoheadrightarrow G_K)$. 
Moreover, there exists a natural isomorphism 
$\Delta _{x_j}\isom \hat {\Bbb Z}(1)$ (cf. loc. cit.).
The profinite group $D_{x_j}$ is isomorphic to the admissible fundamental group $\pi_1^{\adm}(\Cal X)$
of $\Cal X \defeq \Spec \frac {\Cal O_K [[S,T]]}{(ST-\pi)}$, where $\pi$ is a uniformiser of $\Cal O_K$ (cf. loc. cit.). 
The above exact sequence splits. Indeed, the (admissible) covers $\Cal Y_N\to \Cal X$
defined generically by extracting an N-th root of $S$ with $\Cal Y_N$ normal, for all integers $N\ge 1$, define a splitting of this sequence.
Such a splitting induces a section
$s_{x_j}:G_K\to \pi_1(X_{K},\xi)\isom   \pi_1 (X_{k},\xi')^{\adm}$ 
of $\pi_1(X_{K},\xi)$.
The section $s_{x_j}$ is {\it not point-theoretic}. Indeed, if $s_{x_j}$ arises from a rational point $x\in X(K)$, then $x$ 
specialises in a smooth non-cuspidal point of an irreducible component $X_i$ (cf. the discussion before Lemma 2.2.1) 
which is necessarily adjacent to an irreducible component     $X_{i'}$ passing through $x_j$ (cf. [Hoshi-Mochizuki1],
Corollary 1.15(iv)(c)(d)).
Let $s_i:G_k\to \pi_1(U_i,*)^{\tame}$ be the section of $\pi_1(U_i,*)^{\tame}$ which is induced by $s_{x_j}$ (cf. Lemma 2.2.1.).
This section (which is non-cuspidal) would then be cuspidal which contradicts Lemma  2.2.1(ii).

\proclaim {Proposition 2.3.1} Let $K$ be a complete discrete valuation field with residue field $k$ of characteristic $0$. 
Assume that $k$ satisfies the conditions (i) and (iv)(a) (cf. Definition 0.2). Then there exists a flat, proper, smooth, geometrically connected 
hyperbolic curve $C$ over $K$, and a section $s:G_K\to \pi_1(C,*)$ of $\pi_1(C,*)$ which is not point-theoretic.
\endproclaim

\demo {Proof} Write $\Cal O_K$ for the valuation ring of $K$.
Using formal patching techniques one can construct a proper, stable, and regular $\Cal O_K$-curve $X$
satisfying the assumptions in 2.2. In particular, 
$C\defeq X_{K}$ is smooth, hyperbolic, geometrically connected, 
and the double points of $X_k$ are $k$-rational (compare with [Sa\"\i di1], Lemma 6.3). 
As in the above discussion (before Proposition 2.3.1), let $s_{x_j}:G_K\to \pi_1(C,*)$ be a section of $\pi_1(C,*)$
arising from a double point $x_j$ of the special fibre 
$X_s$ of $X$. Then $s_{x_j}$ is not point-theoretic as explained above. 
\qed
\enddemo

\subhead
\S 3.  Sections of arithmetic fundamental groups of curves over function fields in characteristic $0$ 
\endsubhead
We use the notations introduced in $\S1$ and $\S 2$. 
In this section $k$ is a field of characteristic $0$ which {\it strongly satisfies} the condition (i) in Definition 0.2.
Let $C$ be a separated, smooth, and connected 
algebraic curve over $k$ with function field $K\defeq k(C)$.
Let $\Cal X\to C$ be a {\it flat, proper, and smooth relative curve} with generic fibre $X\defeq \Cal X\times _C\Spec K$
which is geometrically 
connected and {\it hyperbolic} (i.e. $g(X)\ge 2$). 
Let $c\in C^{\cl}$ be a closed point and $\Cal X_c\defeq \Cal X\times _C \Spec {k(c)}$ the fibre of $\Cal X$ at $c$. 
Let $\xi$ (resp. $\xi_c$) be a geometric point of $\Cal X$ with value in its generic point (resp. with value in the generic point of $\Cal X_c$).
Then $\xi$ (resp. $\xi_c$) determines an algebraic closure $\overline K$ (resp. $\overline {k(c)}$) of $K$ (resp. of the residue field $k(c)$ of $C$ at $c$).

\proclaim {Lemma 3.1} For $c\in C^{\cl}$ there exists a commutative diagram
$$
\CD
1@>>> \pi_1(X_{\overline K},\overline \xi) @>>> \pi_1(X,\xi) @>>> G_K @>>> 1\\
@.     @V{\id}VV            @VVV           @VVV  @.  \\
1  @>>>   \pi_1(X_{\overline K},\overline \xi) @>>> \pi_1(\Cal X,\xi) @>>> \pi _1(C,\xi) @>>> 1\\
@.     @AAA             @AAA           @AAA  @.  \\
1  @>>>   \pi_1(\Cal X_{\overline {k(c)}},\overline {\xi_c}) @>>> \pi_1(\Cal X_c,\xi_c) @>>> G_{k(c)}@>>> 1\\
\endCD
\tag 3.1$$
where $G_K=\Gal (\overline K/K)$, $G_{k(c)}=\Gal (\overline {k(c)}/k(c))$, $\overline \xi$ (resp. $\overline {\xi_c}$)
are geometric points induced by $\xi$ (resp. $\xi_c$), the middle and right upper vertical maps are natural continuous 
surjective homomorphisms, the middle lower vertical map is defined up to conjugation,
the right lower vertical map is injective, and both squares on the right are cartesian.
 \endproclaim

\demo{Proof} Follows from the functoriality of fundamental groups and the homotopy exact sequence for $\pi_1$ (cf. [Grothendieck1], \'Expos\'e XIII, \S4).
\qed
\enddemo

Let 
$s:G_K\to \pi_1(X,\xi)$ 
be a {\it section} of $\pi_1(X,\xi)$.

\proclaim {Lemma 3.2} There exists a section $s_C:\pi _1(C,\xi)\to  \pi_1(\Cal X,\xi)$ of the projection
 $\pi_1(\Cal X,\xi)\twoheadrightarrow \pi _1(C,\xi)$, which restricts for each closed point $c\in C^{\cl}$ 
 to a section $s_c:G_{k(c)}\to \pi_1(\Cal X_c,\xi_c)$ of $\pi_1(\Cal X_c,\xi_c)$,
and we have a commutative diagram
 $$
 \CD 
G_K @>s>>  \pi_1(X,\xi) \\
@VVV           @VVV \\
\pi _1(C,\xi) @>{s_C}>> \pi_1(\Cal X,\xi) \\
 @AAA           @AAA \\
 G_{k(c)} @>{s_c}>> \pi_1(\Cal X_c,\xi_c) \\
\endCD
\tag {3.2}$$
where the vertical maps are the ones in diagram (3.1).
\endproclaim

\demo{Proof} Follows easily from Lemma 2.1.2 and Lemma 3.1. 
\qed
\enddemo

{\it From now on we assume that $X(K)\neq \emptyset$}. 
Let $\Cal J\defeq \Pic^0_{\Cal X/C}\to C$ be the (relative) jacobian of the (relative) curve $\Cal X\to C$, 
$J\defeq \Cal J_K$ the jacobian variety of $X$,
$\Cal J_c\defeq \Cal J_{k(c)}$
the jacobian variety of $\Cal X_{c}$, and $J_c\defeq J_{K_c}$ the jacobian variety of 
$X_{K_c}$,  where $c\in C^{\cl}$ and $K_c$ is the {\it completion} of $K$ at $c$. The commutative diagram (3.1) induces a commutative diagram of exact sequences of 
{\it geometrically abelian fundamental groups} (cf. $\S1$)
$$
\CD
1@>>> \pi_1(X_{\overline K},\overline \xi)^{\ab} @>>> \pi_1(X,\xi)^{(\ab)} @>>> G_K @>>> 1\\
@.     @V{\id}VV            @VVV           @VVV  @.  \\
1  @>>>   \pi_1(X_{\overline K},\overline \xi)^{\ab} @>>> \pi_1(\Cal X,\xi)^{(\ab)} @>>> \pi _1(C,\xi) @>>> 1\\
@.     @AAA             @AAA           @AAA  @.  \\
1  @>>>   \pi_1(\Cal X_{\overline {k(c)}},\overline {\xi_c})^{\ab} @>>> \pi_1(\Cal X_c,\xi_c)^{(\ab)} @>>> G_{k(c)}@>>> 1\\
\endCD
\tag 3.3$$
satisfying similar properties as in Lemma 3.1. Here $\pi_1(\Cal X,\xi)^{(\ab)}$ is defined as the push forward of the group extension
$1\to \pi_1(X_{\overline K},\overline \xi)\to \pi_1(\Cal X,\xi) \to \pi _1(C,\xi) \to 1$ by the characteristic quotient
$\pi_1(X_{\overline K},\overline \xi)\twoheadrightarrow \pi_1(X_{\overline K},\overline \xi)^{\ab}$.
Fix a section $x\in \Cal X(C)=X(K)$ of $\Cal X\to C$, and the embedding $\Cal X\to \Cal J$ mapping $x$ to the zero section.
We identify $\pi_1(X_{\overline K},\overline \xi)^{\ab}$ (resp.  $\pi_1(\Cal X_{\overline {k(c)}},\overline {\xi_c})^{\ab}$) with the Tate module $TJ$ 
(resp. $T\Cal J_c$), and $TJ_c$ with $T{\Cal J}_c$ ($K_c$ is a complete discrete valuation ring of equal characteristic $0$).

\proclaim{Lemma 3.3} Let $c\in C^{\cl}$. We have a commutative diagram of Kummer exact sequences (cf. $\S1$)
$$
\CD
0 @>>> J(K_c)^{\wedge}   @>>> H^1(G_{K_c},TJ_c)  @>>> TH^1(G_{K_c},J_c) @>>> 0\\
@.    @AAA         @AAA       @AAA \\
0 @>>> \Cal J_c(k(c))^{\wedge} @>>> H^1(G_{k(c)},T{\Cal J}_c)  @>>> TH^1(G_{k(c)},\Cal J_c) @>>> 0\\
\endCD
\tag 3.4$$
where the middle vertical map is the inflation map, the left vertical map is an isomorphism, 
the middle and right vertical maps are injective.
\endproclaim

\demo{Proof} Follows from the fact that there exists an isomorphism
$H^1(G_{k(c)},T{\Cal J}_c)\isom  H^1(\Gal (K_c^{\ur}/K_c), TJ_c)$
where $K_c^{\ur}/K_c$
denotes the maximal unramified  subextension of $\overline {K_c}/K_c$,
and the fact that the kernel of the specialisation map $J_c(K_c)\twoheadrightarrow \Cal J_c (k(c))$ is uniquely divisible 
(cf. [Lang-Tate], Proposition 8). (See also the commutative diagram in loc. cit. page 675.) 
\qed
\enddemo

With the notations in Lemma 3.3, we will identify $H^1(G_{k(c)},T\Cal J_c)$ with its image in $H^1(G_{K_c},TJ_c)$.
Next, we fix the base point of the torsor of splittings of the middle horizontal sequence in diagram (3.3) arising from the section $x\in \Cal X(C)$,
and the corresponding base points of the torsor of splittings of the upper and lower horizontal sequences in diagram (3.3). 
The sections $s:G_K\to \pi_1(X,\xi)$, $s_C:\pi _1(C,\xi)\to \pi_1(\Cal X,\xi)$, and
 $s_c:G_{k(c)}\to \pi_1(\Cal X_c,\xi_c)$ (cf. Lemma 3.2), give rise naturally
to sections  $s^{\ab}:G_K\to \pi_1(X,\xi)^{(\ab)}$, $s_C^{\ab}:\pi _1(C,\xi)\to \pi_1(\Cal X,\xi)^{(\ab)}$,
and $s_c^{\ab}:G_{k(c)}\to \pi_1(\Cal X_c,\xi_c)^{(\ab)}$,
of the projections $\pi_1(X,\xi)^{(\ab)}\twoheadrightarrow G_K$,
$\pi_1(\Cal X,\xi)^{(\ab)}\twoheadrightarrow \pi _1(C,\xi)$, and
$\pi_1(\Cal X_c,\xi_c)^{(\ab)}\twoheadrightarrow G_{k(c)}$; respectively, which 
correspond to elements $s^{\ab}\in H^1(G_{K},TJ)$, $s_C^{\ab}\in H^1(\pi _1(C,\xi), TJ)$,
and $s_c ^{\ab}\in H^1(G_{k(c)},T\Cal J_c)$, for $c\in C^{\cl}$.

\proclaim{Lemma 3.4} Let $c\in C^{\cl}$. We have a commutative diagram of homomorphisms

$$
\CD
H^1(\pi _1(C,\xi), TJ)  @>>>    H^1(G_{K},TJ)\\
@VVV      @VVV\\
H^1(G_{k(c)},T\Cal J_c)   @>>>  H^1(G_{K_c},TJ_c)\\
\endCD
$$
where the horizontal maps are injective inflation maps 
and the vertical maps are restriction homomorphisms, in which $s_C ^{\ab}$ maps to  $s^{\ab}$ (resp. 
$s_c^{\ab}$). In particular, the images of $s^{\ab}$ and $s_c^{\ab}$ in $H^1(G_{K_c},TJ_c)$ are equal.
Moreover, assume that the section $s_c$ is point-theoretic, then 
$s_c^{\ab}$ is point-theoretic (cf. Definition 1.1(iii)). In particular, $s_c^{\ab}$ is pro-geometric both as an element of 
$H^1(G_{k(c)},T\Cal J_c)$ and  $H^1(G_{K_c},TJ_c)$.
(cf. Definition 1.1(i)). 
\endproclaim

\demo {Proof}
The first assertion is easily verified. The second assertion  follows from the various definitions (cf. Lemma 1.2), and Lemma 3.3.
\qed
\enddemo

Next, consider the following commutative diagram
$$
\CD
0 @>>> J(K)^{\wedge}   @>>> H^1(G_{K},TJ)  @>>> TH^1(G_{K},J) @>>> 0\\
@. @VVV        @VVV        @VVV \\
0 @>>> \prod _{c} J_c(K_c)^{\wedge}   @>>> \prod _{c} H^1(G_{K_c},TJ_c)  @>>> \prod _{c}TH^1(G_{K_c},J_c) @>>> 0\\
\endCD
\tag 3.5$$
where the horizontal sequences are the Kummer exact sequences, the vertical maps are the diagonal mappings, and the product in the 
lower sequence is over {\it all} closed points $c\in C^{\cl}$. 
The image of $s^{\ab}\in H^1(G_{K},TJ)$ in $H^1(G_{K_c},TJ_c)$ via the above middle vertical map coincides with the element 
$s_c^{\ab} \in H^1(G_{k(c)},T\Cal J_c)  \subset H^1(G_{K_c},TJ_c)$, $\forall c\in C^{\cl}$ (cf. Lemma 3.4). 
Recall (cf. $\S0$)
$$\Sha (\Cal J)=\Sha (J,C)\defeq \Ker \lgroup H^1(G_{K},J) \to \prod _{c\in C^{\cl}}H^1(G_{K_c},J_c)\rgroup.$$ 
Note that the kernel of the map $TH^1(G_{K},J) \to \prod _{c\in C}TH^1(G_{K_c},J_c)$ in diagram (3.5)
is the Tate module $T\Sha(\Cal J)$ of the Shafarevich-Tate group $\Sha (\Cal J)$.
The following is immediate from the various definitions, and Lemma 3.4.

 \proclaim{Lemma 3.5} Suppose that $s_c^{\ab}\in H^1(G_{K_c},TJ_c)$ 
 is pro-geometric (cf. Definition 1.1(i)) $\forall c\in C^{\cl}$, and $T\Sha(\Cal J)=0$. 
Then $s^{\ab}\in J(K)^{\wedge}$ is pro-geometric. 
 \endproclaim

\subhead
\S 4. Proofs of Theorems A and B
\endsubhead
Next, we prove Theorems A and B.

\subhead
4.1. Proof of Theorem A
\endsubhead
Recall the notations introduced in $\S3$ that we will use throughout. 
Let $k$, $C$, $K=k(C)$, $\Cal X\to C$, and $X\to \Spec K$, be as in $\S3$.
We assume that $X(K)\neq \emptyset$, and $k$ {\it strongly satisfies the conditions} (i), (ii), (iv), (v), and (vi), in Definition 0.2. 
Let $\Cal J\defeq \Pic^0_{\Cal X/C}$ be the relative jacobian of $\Cal X\to C$, and assume that $T\Sha(\Cal J)=0$.
We will show that the map (cf. $\S0$)
$$\varphi_X: X(K)\to  \overline {\Sec} _{\pi _1(X,\xi)},\ \ \ \ x\mapsto \varphi_X (x)=[s_x],$$
is {\it bijective}.

First, we prove $\varphi_X$ is {\it injective}. Let $x_1,x_2\in X(K)$
such that $[s_1\defeq s_{x_1}]=[s_2\defeq s_{x_2}]$ holds in $\overline {\Sec} _{\pi _1(X,\xi)}$.
Thus, $[s_{1,c}]=[s_{2,c}]$ in $\overline {\Sec} _{\pi _1(\Cal X_c,\xi_c)}$, $\forall c\in C^{\cl}$ (cf. Lemma 3.2), and $s_{i,c}=s_{x_{i,c}}$ is point-theoretic, 
where $x_{i,c}\in \Cal X_c(k(c))$ is uniquely determined by $s_{i,c}$,
for $i\in \{1,2\}$, as we assumed that $k$ strongly satisfies (ii). 
The map $\varphi_{\Cal X_c}: \Cal X_c(k(c))\to  \overline {\Sec} _{\pi _1(\Cal X_c,\xi_c)}$
is bijective by assumption, hence $x_{1,c}=x_{2,c}$, $\forall c\in C^{\cl}$. Moreover, $x_{i,c}$ is the specialisation of $x_i$ in $\Cal X_c$ (cf. discussion before Lemma 2.1.3). 
From this it follows that $x_1=x_2$
and $\varphi_X$ is injective. Indeed, the natural specialisation map $X(K)\to \prod _{c\in C^{\cl}} \Cal X_c(k(c))$ is injective.

Next, we prove that $\varphi_X$ is {\it surjective}.
Let 
$$s:G_K\to \pi_1(X,\xi)$$ 
be a {\it section} of $\pi_1(X,\xi)$.
We will show that $s$ is {\it point-theoretic} under the above assumptions.
First, we have a diagram (3.2) and it follows from the condition (ii) that the section $s_c$ (cf. loc. cit.) 
is point-theoretic and arises from a unique rational point $x_c\in \Cal X(k(c))$, $\forall c\in C^{\cl}$.
Moreover, $s_c^{\ab}\in J(K_c)^{\wedge}\subset H^1(G_{K_c},TJ_c)$ is pro-geometric in the sense of 
Definition 1.1(i), $\forall c\in C^{\cl}$ (cf. Lemma 3.4). Then it follows from the assumption
$T\Sha(\Cal J)=0$ that $s^{\ab}\in J(K)^{\wedge}\subset H^1(G_{K},TJ)$ is pro-geometric (cf. Lemma 3.5).

\proclaim{Lemma 4.1.1} The natural homomorphism $J(K)\to J(K)^{\wedge}$ is injective and 
$s^{\ab}\in J(K)\subseteq J(K)^{\wedge}$             
is geometric.
\endproclaim

\demo{Proof} There exist closed points $c_1,c_2\in C^{\cl}$ such that the natural specialisation 
homomorphism $J(K)\to \Cal J_{c_1}(k(c_1))\times  \Cal J_{c_2}(k(c_2))$ is injective (cf. [Poonen-Voloch], Proposition 2.4). 
We have a commutative diagram of exact sequences

$$
\CD
0 @>>>  J(K) @>>> H\defeq \Cal J_{c_1}(k(c_1))\times  \Cal J_{c_2}(k(c_2)) @>>>  H/J(K) @>>> 0\\
@.    @VVV @V{\phi} VV @VVV  @. \\
0 @>>> J(K)^{\wedge}  @>{\psi}>>  H^{\wedge}=\Cal J_{c_1}(k(c_1))^{\wedge}\times  \Cal J_{c_2}(k(c_2))^{\wedge}   @>>> (H/J(K)) ^{\wedge} @>>> 0  \\
\endCD
$$
where the right and middle vertical maps are injective by the assumption (iv)(a), 
and the maps $\psi$ and $\phi$ are the natural ones. (The exactness of the lower sequence in the above diagram follows easily from the assumption (iv)(b).)
In particular, the left vertical map is injective, and the equality
$J(K)=\phi(H)\cap \psi(J(K)^{\wedge})$ holds inside $H^{\wedge}$. 
The image of $s^{\ab}\in J(K)^{\wedge}$ in $H^{\wedge}$ via the map $\psi$ is the element 
$(s_{c_1}^{\ab}, s_{c_2}^{\ab})  \in H^{\wedge}\subset H^1(G_{k(c_1)},T\Cal J_{c_1})\times H^1(G_{k(c_2)},T\Cal J_{c_2})$
associated to the sections
$s_{c_i}:G_{k(c_i)}\to \pi_1(\Cal X_{c_i},\xi_{c_i})$ (induced by the section $s$), for $i\in \{1,2\}$, 
which are point-theoretic by the assumption (ii).
In particular, $s_{c_i}^{\ab}\in \Cal J_{c_i}(k(c_i))$ is geometric (cf. Lemma 1.2), and $s^{\ab}\in J(K)$ is geometric by the above discussion.
\qed
\enddemo

We fix an embedding  $\Cal X@>\iota>> \Cal J$ mapping an element of $\Cal X(C)=X(K)$ (which is non-empty by
our assumptions) to the zero section.

\proclaim {Lemma 4.1.2} The element  $s^{\ab}\in \iota(X(K))\subset J(K)$ is point-theoretic. 
\endproclaim

\demo{Proof} For each closed point $c\in C^{\cl}$ the element $s_c^{\ab}\in \Cal J (k(c))^{\wedge}\subset H^1(G_{k(c)},T\Cal J_c)$ 
corresponding to the section $s_c:G_{k(c)}\to \pi_1(\Cal X_c,\xi_c)$
lies in the subset $\iota (\Cal X(k(c)))\subset \Cal J(k(c))\subset \Cal J (k(c))^{\wedge}$, since the section $s_c$ is point-theoretic by (ii). 
We view $s^{\ab}\in J(K)$ as a rational section of $\Cal J\to C$; in fact $s^{\ab}:C\to \Cal J$ is a morphism since $C$ is a smooth curve. 
For each closed point
$c\in C^{\cl}$ the image $s^{\ab}(c)$ is a closed point of $\Cal X_c\subset \Cal J_c$, where we view $\Cal X_c$ as a closed subscheme of $\Cal J_c$ 
via the closed immersion $\Cal X_{c} @>{\iota _c}>> \Cal J_c$ induced by $\iota$. From this it follows that the morphism $s^{\ab}: C\to \Cal J$
factorizes as $s^{\ab}: C\to \Cal X @>\iota>> \Cal J$ and $s^{\ab}$ belongs to the subset $\iota (X(K))\subseteq J(K)$. 
\qed
\enddemo

Let $\tilde x\in X(K)$ such that $\iota (\tilde x)=s^{\ab}$ (cf. Lemma 4.1.2).
For $c\in C^{\cl}$ let $\tilde x_c$ be the specialisation of $\tilde x$ in $\Cal X_c$.

\proclaim {Lemma 4.1.3} The equality $\tilde x_c=x_c$ holds in $\Cal X(k(c))$, $\forall c\in C^{\cl}$.
\endproclaim

\demo{Proof} First, the equality $s_c^{\ab}=s_{\tilde x_c}^{\ab}=s_{x_c}^{\ab}$ holds in $H^1(G_{k(c)},T\Cal J_c)$ (cf. Lemma 3.4).
The equality $\tilde x_c=x_c$ follows then from the injectivity of the maps $\iota (\Cal X_{c})\hookrightarrow \Cal J_c(k(c))
\hookrightarrow \Cal J_c(k(c))^{\wedge}\to H^1(G_{k(c)},T\Cal J_c)$, for $c\in C^{\cl}$ (see the assumption (iv)(a)). 
\qed
\enddemo

Next, and in order to show that the section $s$ is point-theoretic it suffices to show, by a well-know {\it limit argument} in anabelian geometry 
(cf. [Tamagawa], Proposition 2.8 (iv)), using the assumption (v) (cf. Definition 0.2), the following.
Let $H\subseteq \pi_1(X,\xi)$ be an open subgroup such that $s(G_K)\subset H$, corresponding to an \'etale cover $Y\to X$, 
then $Y(K)\neq \emptyset$ holds.
There is a natural identification $H=\pi_1(Y,\xi)$. Moreover, the cover $Y\to X$ extends to an \'etale cover
$\Cal Y\to \Cal X$ and $\Cal Y$ is a smooth model of $Y$ over $C$ (cf. Lemma 3.1 and Lemma 3.2). 
For a closed point $c\in C^{\cl}$ write $\Cal Y_c\defeq \Cal Y_{k(c)}$.
We will show that $Y(K)\neq \emptyset$.

Let $s':G_K\to  \pi_1(Y,\xi)$ be the section of $\pi_1(Y,\xi)=H$ induced by $s$,
which extends to a section $s':\pi_1(C,\xi)\to  \pi_1(\Cal Y,\xi)$ of the projection $\pi_1(\Cal Y,\xi)\twoheadrightarrow \pi_1(C,\xi)$, and further induces
a section $s'_c:G_{k(c)}\to  \pi_1(\Cal Y_c,\xi_c)$ of $\pi_1(\Cal Y_c,\xi_c)$, $\forall c\in C^{\cl}$ (cf. loc. cit.). Note that 
$s'_c$ is induced by the section $s_c$. The section $s'_c$ is point-theoretic and arises from a unique rational point 
$y_c\in \Cal Y_c(k(c))$. Moreover, $x_c$ is the image of $y_c$ in $\Cal X_c$
via the morphism $\Cal Y_c\to \Cal X_c$ (cf. condition (ii), the fact that $s_c$ is point-theoretic and arises from $x_c$, and $s'_c$ is induced by $s_c$).
View $\tilde x\in X(K)=\Cal X(C)$ as a section $\tilde x:C\to \Cal X$, and let $\Cal Y_{\tilde x}$ be the scheme-theoretic inverse image of $\tilde x(C)$ in $\Cal Y$ 
via the above \'etale map 
$\Cal Y\to \Cal X$. Thus, $\Cal Y_{\tilde x}\to \tilde x(C)$ is a finite \'etale map. We have $y_c\in \Cal Y_{\tilde x}(k(c))$, $\forall c\in C^{\cl}$, as follows from the various definitions. 
Then $\Cal Y_{\tilde x}(K)\neq \emptyset$ by the assumption (vi), and a fortiori $\Cal Y_{\tilde x}(K)\subseteq \Cal Y(K)=Y(K)\neq \emptyset$.

Thus, we proved that $[s]=[s_x]$ holds in $\overline {\Sec} _{\pi _1(X,\xi)}$ for a (unique) $x\in X(K)$. The following follows from Lemma 4.1.3 (cf. above proof that $\varphi_X$
is injective).

\proclaim {Lemma 4.1.4} The equality $x=\tilde x$ holds.
\endproclaim

This finishes the proof of Theorem A.
\qed

Finally, we show that Hilbertian fields satisfy the condition (vi).

\proclaim{Lemma 4.1.5}
Let $k$ be a Hilbertian field. Then $k$ strongly satisfies the condition (vi).
\endproclaim

\demo{Proof} Let $k'/k$ be a finite extension and $C$ a separated, smooth, and connected curve over $k'$ with function field $K\defeq k'(C)$.
Let $\widetilde C\to C$ be a finite \'etale cover with $\widetilde C_c(k'(c))\neq \emptyset$, $\forall c\in C^{\cl}$ with residue field $k'(c)$.
Note that $k'$ is Hilbertian (cf. [Serre], 9.5). We show $\widetilde C(K)\neq \emptyset$.
Assume that $\widetilde C(K)=\emptyset$. Then for each connected component $\widetilde C_{\alpha}$ of $\widetilde C$ the degree of the morphism 
$\widetilde C_{\alpha}\to C$ is $\ge2$. Hilbert's irreducibility Theorem (cf. [Serre], 9.2) implies that there exist points $c\in C^{\cl}$ such that the fibre of 
$c$ in each connected component of $\widetilde C$ is irreducible. This contradicts the assumption that $\widetilde C_c(k(c))\neq \emptyset$, $\forall c\in C^{\cl}$. 
Thus, $\widetilde C(K)\neq \emptyset$ holds.
\qed
\enddemo

In the course of proving Theorem A we proved the following \lq\lq{\it adelic Mordell-Lang}\rq\rq statement.
\proclaim{Proposition 4.1.6} With the same notations as above, 
assume that $k$ only satisfies the condition (iv) in Definition 0.2 (where we take $k'=k$), and $\Cal X(C)\neq \emptyset$.
Then the map $J(K)^{\wedge}\to  \prod _{c\in C^{\cl}} J_c(K_c)^{\wedge}\isom  \prod _{c\in C^{\cl}} \Cal J_c(k(c))^{\wedge}$ 
(cf. diagram (3.4) and Lemma 3.3) is injective. Further, inside $\prod _{c\in C^{\cl}} \Cal J_c(k(c))^{\wedge}$
the equality $J(K)^{\wedge}\bigcap \lgroup \prod _{c\in C^{\cl}} \iota (\Cal X_c(k(c)))\rgroup=\iota (X(K))$ 
holds.
\endproclaim

\demo{Proof} See the statements and proofs of lemma 4.1.1 and Lemma 4.1.2.
\qed
\enddemo

\subhead
4.2. Proof of Theorem B
\endsubhead
We briefly explain how Theorem B can be derived from Theorem A.
Let $k$ be a field with $\char (k)=0$,  and $K=k(C)$ the function field of a separated, smooth, and connected curve $C$ over $k$.
Assume that $k$ satisfies the condition $(\star)$ (cf. Definition 0.2). Let $L/K$ be a finite extension. Then we prove the {\bf SC} holds over $L$.
We can, without loss of generality, assume that $L=K$. Let $X$ be a proper, smooth, and hyperbolic curve over $K$, we need to prove that the {\bf SC}
holds for $X$, i.e., that the map (cf. $\S0$)
$\varphi_X: X(K)\to  \overline {\Sec} _{\pi _1(X,\xi)}$
is {\it bijective}.
The injectivity of the map $\varphi_X$ follows as in the proof of Theorem A, where one only uses the fact that $k$ strongly satisfies the condition (ii) (cf. loc. cit.)

Next, we prove that $\varphi_X$ is {\it surjective}. Let 
$s:G_K\to \pi_1(X,\xi)$ 
be a {\it section} of $\pi_1(X,\xi)$.
We need to show that $s$ is {\it point-theoretic} under the above assumptions.
For this purpose we can, in the course of the proof,  replace $K$ by a finite extension $L'/K$.
Indeed, let $L'/K$ be a finite Galois extension,
$s':G_{L'}\defeq \Gal (\overline K/L')\to \pi_1(X_{L'},\xi)$ the section of $\pi_1(X_{L'},\xi)$
which is induced by $s$, and assume that $s'=s_{x'}$ is point-theoretic where $x'\in X(L')$. Then $s'(G_{L'})$ is self-normalising
in  $\pi_1(X_{L'},\xi)$, and $s(G_K)$ is contained in the normaliser of $s'(G_{L'})$ in  $\pi_1(X,\xi)$
which coincides with a decomposition group associated to the image $x$ of $x'$ in $X$
(this follows from condition $(\star)(iv)(a)$, cf. proof of Lemma 2.2.1 and the references therein). 
The point $x$ is then necessarily $K$-rational. 
Now consider a finite (Galois) extension $L'/K$ such that $X(L')\neq \emptyset$, then the corresponding section
$s':G_{L'}\defeq \Gal (\overline K/L')\to \pi_1(X_{L'},\xi)$ of $\pi_1(X_{L'},\xi)$ which is induced by $s$ is point-theoretic by Theorem A
and our assumption that $k$ satisfies the condition $(\star)$.
From this it follows that $s$ is point theoretic (cf. above discussion).
\qed

\subhead
\S 5.  Sections of arithmetic fundamental groups of curves over finitely generated fields in characteristic $0$ 
\endsubhead
In this section we combine the techniques/results in $\S3$ and $\S4$ with the results in [Sa\"\i di-Tamagawa] to prove variants of 
Theorems A and B, as well as Theorem C.
We use the notations introduced in previous sections. All fields in this section are {\it finitely generated over $\Bbb Q$}.
In order to deduce the validity of the {\bf SC} over finitely generated fields from its validity over number fields (which for the time being is still not known), 
one is reduced (by an induction argument)
to prove the following.

\proclaim {Conjecture D} Assume that the {\bf SC} holds over all finitely generated fields of transcendence degree $i\ge 0$. then the {\bf SC} 
holds over all finitely generated fields of transcendence degree $i+1$.
\endproclaim

Next, we recall the {\it discrete Selmer conjecture} formulated in [Sa\"\i di-Tamagawa]. Let $k$ be a {\it finitely generated} field over $\Bbb Q$, $C$ a separated, smooth, and geometrically connected curve over $k$, $K\defeq k(C)$ the function field of $C$, and 
$$\Cal A\to C$$ 
an {\it abelian scheme} over $C$. 
Write $A\defeq \Cal A\times _C \Spec K$ for the generic fibre of $\Cal A$.
For a {\it closed} point $c\in C^{\cl}$, let $k(c)$ be the residue field at $c$, and $\Cal A_c\defeq \Cal A_{k(c)}=\Cal A\times _C\Spec k(c)$.
Then we have a commutative diagram of exact sequences

$$
\CD
0 @>>> A(K)^{\wedge} @>>> H^1(G_K,TA) @>>> TH^1(G_K,A) @>>> 0  \\
@. @VVV      @V\prod_c \res_cVV    @V\prod_c\res_cVV \\
0 @>>> \prod_c\Cal A_c(k(c))^{\wedge} @>>>  \prod_c H^1(G_{k(c)},T \Cal A_c) @>>> \prod_c TH^1(G_{k(c)},\Cal A_c) @>>> 0 \\
\endCD
$$
where the horizontal sequences are the Kummer exact sequences, the vertical maps are natural 
restriction homomorphisms, and the product is taken over all closed points $c\in C^{\cl}$ (cf. loc. cit. diagram (0.1)).
Define the {\it profinite Selmer group} 
$$\Sel (\Cal A)\defeq \Sel (A,C)\defeq \Ker \lgroup H^1(G_K,TA) \to \prod_c TH^1(G_{k(c)},\Cal A_c)\rgroup,$$
and the {\it Shafarevich-Tate group}
$$\Sha(\Cal A)\defeq \Sha (A,C)\defeq \Ker \lgroup H^1(G_K,A)\to \prod_c H^1(G_{k(c)},\Cal A_c)\rgroup,$$
where the product is taken over {\it all} closed points $c\in C^{\cl}$.
For each $c\in C^{\cl}$ the group $\Cal A(k(c))$ of $k(c)$-rational points of $\Cal A$
is {\it finitely generated} as $k(c)$ is finitely generated over $\Bbb Q$ (Mordell-Weil Theorem), hence injects into its profinite completion
$\Cal A(k(c))^{\wedge}$. We identify $\Cal A(k(c))$ with its image in $\Cal A(k(c))^{\wedge}$.
Define the {\it discrete Selmer group} by
$$\Se (\Cal A)\defeq \Se (A,C)\defeq \Sel (\Cal A) \bigcap \prod_c \Cal A(k(c)) \subset \prod_cH^1(G_{k(c)},T\Cal A_c),$$
where the product is taken over {\it all} closed points $c\in C^{\cl}$.
Note that $A(K)\subseteq \Se(\Cal A)$. In [Sa\"\i di-Tamagawa] we conjectured the following (cf. loc. cit. Conjecture E).

\definition {Conjecture E}
The equality $\Se(\Cal A)=A(K)$ holds.
\enddefinition

Moreover, we proved the following facts.
\definition {Facts}
\noindent
{\bf (i)}\  $\Se (\Cal A)$ is a {\it finitely generated} abelian group (cf. loc. cit. Proposition 2.5).

\noindent
{\bf (ii)}\ Assume that there exists a prime integer $\ell$ such that the $\ell$-primary part 
$\Sha (\Cal A) [\ell ^{\infty}]\defeq \bigcup _{n\ge 1}\Sha (\Cal A)[l^n]$ 
of $\Sha(\Cal A)\defeq \Sha (A,C)$ is {\it finite}. Then 
$\Se(\Cal A)=A(K)$, i.e., Conjecture E holds in this case (cf. loc. cit. Proposition 3.7).

\noindent
{\bf (iii)}\ Assume that $A$ is {\it isotrivial}, i.e., $A_{\overline K}$ is defined over $\bar k$, then  
$\Sha(\Cal A)$ is {\it finite} (cf. loc. cit. Theorem 4.2). In particular,
$\Se(\Cal A)=A(K)$ holds, i.e., Conjecture E holds in this case (cf. above fact (ii)).
\enddefinition

One of our main results in this paper is a (conditional) proof of Conjecture D. More precisely, 
we prove the following (cf. discussion after Corollary 5.3).

\proclaim {Theorem 5.1} Assume that Conjecture {\bf E} holds. Then Conjecture {\bf D} holds true.
\endproclaim

In other words, Theorem 5.1 asserts that the validity of the above discrete Selmer conjecture implies that one can reduce the section conjecture {\bf SC} over {\it all finitely generated fields} to the case of {\it number fields}.
The following Corollaries follow immediately from Theorem 5.1, and the above cited Facts (ii), and (iii); respectively. (For the proof of Corollary 5.3, 
use an induction argument.)

\proclaim {Corollary 5.2} Assume that the following holds. Given a finitely generated field $k$ over $\Bbb Q$, a separated, smooth, and connected curve
$C$ over $k$ with generic point $\eta$, and an abelian scheme $\Cal A\to C$ with generic fibre $A\defeq \Cal A\times_C \eta$, then 
there exists a prime integer $\ell$ such that the $\ell$-primary part $\Sha (\Cal A) [\ell ^{\infty}]\defeq \bigcup _{n\ge 1}\Sha (\Cal A)[l^n]$ of $\Sha(\Cal A)\defeq 
\Sha(A,C)$ is {\it finite}.  Then Conjecture {\bf D} holds true.
\endproclaim

\proclaim {Corollary 5.3} Let $K$ be a finitely generated field with algebraic closure $\overline K$ and $\overline {\Bbb Q}$ the algebraic closure of $\Bbb Q$ in $\overline K$. Assume that the {\bf SC} holds over all number fields. Then the {\bf SC} holds for every projective smooth and geometrically connected hyperbolic curve $X$ over $K$ 
such that $X\times _{\Spec K} {\Spec \overline K}$ is defined over $\overline {\Bbb Q}$. 
\endproclaim

The rest of this section is devoted to the proof of Theorem 5.1. To this end it suffices to prove the following (compare with Theorem A). 
(Note that the assumption $X(K)\neq \varnothing$ in Theorem 5.4 is not restrictive in order to deduce Theorem 5.1 
from Theorem 5.4 (cf. arguments in the proof of Theorem B).)

\proclaim {Theorem 5.4} Assume that the {\bf SC} holds over all finitely generated fields of transcendence degree $i\ge 0$. Let $k$ be a field with $\tr\deg_{\Bbb Q}k=i$,
and $K=k(C)$ the function field of a separated, smooth, and connected curve $C$ over $k$. Let $\Cal X\to C$ be a flat, proper, and smooth relative curve with generic fibre
$X\defeq \Cal X\times _C\Spec K$ which is hyperbolic and
geometrically connected, $J\defeq \Pic^0_{X/K}$, and $\Cal J\defeq \Pic^0_{\Cal X/C}$. Assume that $X(K)\neq \emptyset$, and either there exists a prime integer $\ell$ such that the 
$\ell$-primary part $\Sha (\Cal J) [\ell ^{\infty}]$ of $\Sha(\Cal J)\defeq \Sha (J,C)$ is finite or $\Se(\Cal J)\defeq \Se(J,C)=J(K)$ holds. Then the {\bf SC} holds for $X$ (over $K$).
\endproclaim

\demo {Proof} The proof of the injectivity of the map 
$\varphi_X: X(K)\to  \overline {\Sec} _{\pi _1(X,\xi)}$ is similar to the proof of injectivity in the statement of Theorem A (cf. discussion before Lemma 4.1.1).
In what follows we prove $\varphi$ is surjective. Let 
$s:G_K\to \pi_1(X,\xi)$ be a {\it section} of $\pi_1(X,\xi)$. We will show that $s$ is {\it point-theoretic} under the above assumptions.
We use the notations in $\S3$.
Recall that $s$ induces a section $s_c:G_{k(c)}\to \pi_1(\Cal X_c,\xi_c)$ of $\pi_1(\Cal X_c,\xi_c)$, $\forall c\in C^{\cl}$ (cf. Lemma 3.2).
By our assumption that the {\bf SC} holds over all finitely generated fields of transcendence degree $i$ we deduce that the section $s_c$ is {\it point-theoretic}
and arises from a unique rational point $x_c\in \Cal X(k(c))$, $\forall c\in C^{\cl}$. The element $s_c^{\ab}\in H^1(G_{k(c)}, T\Cal J_c)$ is then point-theoretic, hence is geometric
(cf. Lemma 1.2). Thus, $s_c^{\ab}\in \Cal J(k(c))\subseteq \Cal J(k(c))^{\wedge}\subseteq H^1(G_{k(c)}, T\Cal J_c)$, $\forall c\in C^{\cl}$.
In particular, $s^{\ab}\in \Se(\Cal J)\defeq \Se(J,C)$ is in the {\it discrete Selmer group}.
By our second assumption, that either there exists a prime integer $\ell$ such that the 
$\ell$-primary part $\Sha (\Cal J) [\ell ^{\infty}]$ of $\Sha(\Cal J)\defeq \Sha (J,C)$ is finite or $\Se(\Cal J)\defeq \Se(J,C)=J(K)$ holds (note that the former condition implies the latter one by Fact (ii) above), we deduce that $s^{\ab}\in J(K)$ is geometric.
Moreover, a similar proof to that of Lemma 4.1.2 shows that $s^{\ab}$ lies in the subset $X(K)$ of $J(K)$ (cf. loc. cit.). The rest of the proof is similar to that of Theorem A 
(cf. Lemma 4.1.3 and Lemma 4.1.4).

This finishes the proof of Theorem 5.4, and Theorem 5.1.  
\qed
\enddemo

\definition {Remark 5.5} The condition in Theorem 5.4  involving $\Sha (\Cal J) [\ell ^{\infty}]$ is weaker than the condition in
Theorem A and Corollary A involving $\Sha (\Cal J)$.
\enddefinition

%\bigskip
%\bigskip
%\newpage
$$\text{References.}$$

\noindent
[Cadoret-Tamagawa] Torsion of abelian schemes and rational points on moduli spaces, In: Algebraic Number Theory and Related Topics 2007,
Res. Inst. Math. Sci. (RIMS) K\^oky\^uroku Bessatsu B12, Kyoto, 7-29, 2009. 

%\noindent
%[Fresnel-Van der Put] Fresnel, J., Van der Put, M., Rigid analytic geometry and its applications, Progress in Mathematics 218, Birkh\"auser, 2004. 

\noindent
[Grothendieck] Grothendieck, A., Brief an G. Faltings, (German), with an
english translation on pp. 285-293,
London Math. Soc. Lecture Note Ser., 242, Geometric Galois actions, 1,
49-58, Cambridge Univ. Press,
Cambridge, 1997.

\noindent
[Grothendieck1] Grothendieck, A., Rev\^etements \'etales et groupe fondamental, Lecture 
Notes in Math. 224, Springer, Heidelberg, 1971.

\noindent
[Hoshi-Mochizuki] Hoshi, Y., Mochizuki, S., On the combinatorial anabelian geometry of nodally nondegenerate outer representations, Hiroshima Math. J. 41 (2011), no.3, 375-342.

\noindent
[Hoshi-Mochizuki1] Hoshi, Y., Mochizuki, S., Topics surrounding the combinatorial anabelian geometry of hyperbolic curves IV: discretness and sections. Manuscript, available
in http://www.kurims.kyoto-u.ac.jp/~motizuki/papers-english.html, (2015).

\noindent
[Lang-Tate] Lang, S., Tate, J., Principal homogeneous spaces over abelian varieties, American Journal of Mathematics,
Volume 80, Issue 3, (1958), 659-684.

\noindent
[Liu] Liu, Q., Algebraic geometry and arithmetic curves, Oxford graduate texts in mathematics 6. Oxford University Press, 2002.

\noindent
[Mochizuki] Mochizuki, S., The local pro-$p$ anabelian geometry of curves, Invent. Math.  138  (1999), 
no. 2, 319--423.

\noindent
[Mochizuki1] Mochizuki, S., The profinite Grothendieck conjecture for closed hyperbolic curves over number fields, 
J. Math. Sci., Univ. Tokyo 3 (1996), 571-627.

\noindent
[Mochizuki2] Mochizuki, S., Galois sections in absolute anabelian geometry, Nagoya Math. J. 179 (2005), pp. 17-45.

\noindent
[Poonen-Voloch] Poonen, B., Voloch, P., The Brauer-Manin obstruction for subvarieties of abelian varieties over function fields,
Ann. of Math. (2) 171 (2010), no. 1, 511-532.

\noindent
[Sa\"\i di] Sa\"\i di, M., The cuspidalisation of sections of arithmetic fundamental groups, Advances in Mathematics, 230 (2012) 1931-1954.

\noindent
[Sa\"\i di1] Sa\"\i di, M., Rev\^etements \'etale ab\'eliens, courants sur les graphes et 
r\'eduction semi-stable des courbes, Manuscripta Math. 89 (1996), no.2, 245-265.

\noindent
[Sa\"\i di2] Sa\"\i di, M., On the Grothendieck anabelian section conjecture over finitely generated fields. 
Talk at the Utah Summer Research Institute on Algebraic Geometry 2015, video talk available on http://www.claymath.org/utah-videos.

\noindent
[Sa\"\i di-Tamagawa] Sa\"\i di, M., Tamagawa, A., On the arithmetic of abelian varieties, arxiv:1512.00773.

\noindent
[Serre] Serre, J-P., Lectures on the Mordell-Weil Theorem, Translated and edited by Martin Brown from notes by Michel Waldschmidt, $2^{\nd}$ edition, 
Aspect of Mathematics, 1990.

\noindent
[Tamagawa] Tamagawa, A., The Grothendieck conjecture for affine curves,  Compositio Math.  109  (1997),  no. 2, 135--194.

\bigskip

\noindent
Mohamed Sa\"\i di

\noindent
College of Engineering, Mathematics, and Physical Sciences

\noindent
University of Exeter

\noindent
Harrison Building

\noindent
North Park Road

\noindent
EXETER EX4 4QF 

\noindent
United Kingdom

\noindent
M.Saidi\@exeter.ac.uk

\enddocument
\end